\documentclass[11pt]{amsart}

\usepackage{hyperref,amsmath,amssymb,amsthm,tikz}

\newtheorem{thm}{Theorem}
\newtheorem{prop}[thm]{Proposition}
\newtheorem{lem}[thm]{Lemma}
\newtheorem{cor}[thm]{Corollary}

\newtheorem{rem}[thm]{Remark}
\newtheorem{conj}[thm]{Conjecture}

\title{Bounds on the closeness centrality of a graph}

\author{%
Thomas Britz, 
Xin Hu, 
Abdellah Islam and 
Hopein Christofen Tang 
}

\date{}

\begin{document}

\begin{abstract}
We present new values and bounds on the (normalised) closeness centrality $\bar{\mathsf{C}}_C$ of connected graphs
and on its product $\bar{l}\bar{\mathsf{C}}_C$ with the mean distance $\bar{l}$ of these graphs.
Our main result presents the fundamental bounds $1\leq \bar{l}\bar{\mathsf{C}}_C<2$.
The lower bound is tight and the upper bound is asymptotically tight.
Combining the lower bound with known upper bounds on the mean distance,
we find ten new lower bounds for the closeness centrality of graphs.
We also present explicit expressions for $\bar{\mathsf{C}}_C$ and $\bar{l}\bar{\mathsf{C}}_C$
for specific families of graphs. 
Elegantly and perhaps surprisingly, 
the asymptotic values 
$n\bar{\mathsf{C}}_C\big(P_n\big)$ and of 
$n\bar{\mathsf{C}}_C\big(L_n\big)$ 
both equal~$\pi$, 
and the asymptotic limits of $\bar{l}\bar{\mathsf{C}}_C$ 
for these families of graphs are both equal to~$\pi/3$.
We conjecture that the set of values $\bar{l}\bar{\mathsf{C}}_C$ for all connected graphs is dense in the interval $[1,2)$.
\end{abstract}

\maketitle 


\section{Introduction}
\label{sec:intro}
 
\noindent
Let $G = (V,E)$ be a connected graph on $|V| = n$ vertices.
Defined in 1950 by Bavelas~\cite{bavelas50},
the {\em (normalised) closeness centrality} of~$G$ 
is defined as the average of the (normalised) closeness centralities $\bar{\mathsf{C}}_C(v)$ for each vertex $v$ in~$G$:
\[
     \bar{\mathsf{C}}_C(G)
  := \dfrac{1}{n}\sum_{v\in V}\bar{\mathsf{C}}_C(v) \qquad\text{where}\qquad
     \bar{\mathsf{C}}_C(v)
   = \dfrac{n-1}{\displaystyle\sum_{w\in V}d(v,w)}
\]
and where $d(v,w)$ is the {\em distance} between $v$ and~$w$; 
that is, the shortest length of a path between $v$ and~$w$ in~$G$.

\medskip

The purpose of this paper is to provide values and bounds on $\bar{\mathsf{C}}_C:=\bar{\mathsf{C}}_C(G)$
and on the product $\bar{l}\bar{\mathsf{C}}_C$ 
where 
$\bar{l}:=\bar{l}(G)$ is the {\em mean distance} of~$G$, defined as
\[
     \bar{l}(G)
   = \dfrac{1}{n(n-1)}\sum_{u,v\in V} d(u,v)\,.
\]
In Section~\ref{sec:duality}, 
Theorem~\ref{thm:CC-l-inequalities} presents the fundamental and elegant new bounds
$1\leq\bar{l}\bar{\mathsf{C}}_C< 2$.
The lower bound is tight and the upper bound is asymptotically tight.
Although these bounds are simple, they appear to have been overlooked in previous literature.
They show that the closeness centrality $\bar{\mathsf{C}}_C(G)$ and the mean distance $\bar{l}(G)$ 
cannot both be small nor both be large.
These bounds are also useful: 
the lower bound combined with known upper bounds on $\bar{l}(G)$, 
the yields ten new lower bounds for $\bar{\mathsf{C}}_C(G)$;
see Corollary~\ref{cor:lower bounds on CC}.

In Section~\ref{sec:upper-bounds},
we prove four new upper bounds on $\bar{\mathsf{C}}_C(G)$
that depend on different graph parameters;
see Theorems~\ref{thm:edge-deletion}, \ref{thm:non-trivial upper bound 1} and~\ref{thm:non-trivial upper bound 2}.

In Section~\ref{sec:families}, 
we provide lower bounds for three infinite families of graphs, 
namely the self-complementary graphs, the 2-connected graphs, and trees;
see Theorem~\ref{thm:lower bounds on CC for specific graphs}.
We also provide explicit expressions for the closeness centralities $\bar{\mathsf{C}}_C(v)$
for the vertices in twelve specific infinite families of graphs, 
as well as the closeness centralities $\bar{\mathsf{C}}_C(G)$ for ten of these graphs; 
see Proposition~\ref{prop:closeness-centrality-for-various-graph-families}
and Corollary~\ref{cor:closeness-centrality-for-families-of-graphs}.
For the remaining two families, namely the path graphs $P_n$ and the ladder graphs~$L_n$, 
we provide asymptotically tight lower and upper bounds; 
see Propositions~\ref{prop:cool-lower-and-upper-bounds-on-CC-path} 
and~\ref{prop:cool-lower-and-upper-bounds-on-CC-ladder}.
From these bounds, we deduce Corollary~\ref{cor:CCPn-CCLn}
which states the elegant and perhaps surprising fact that the asymptotic limits of 
$n\bar{\mathsf{C}}_C\big(P_n\big)$ and of 
$n\bar{\mathsf{C}}_C\big(L_n\big)$ 
both equal~$\pi$.

In Section~\ref{sec:product-of-closeness-and-mean-distance},
we derive the precise values of the product $\bar{l}\bar{\mathsf{C}}_C$
for the twelve families of graphs in Section~\ref{sec:families} 
as well as another family of graphs.
Again, we find surprising asymptotic values, 
including the value $\displaystyle\lim_{n\to\infty}\bar{l}\bar{\mathsf{C}}_C=\frac{\pi}{3}$ for path and ladder graphs;
see Corollary~\ref{cor:CC-l-converge to pi/3}.

The final Section~\ref{sec:proofofmaintheorem} is dedicated to proving Theorem~\ref{thm:CC-l-inequalities}.
To do so, 
we first prove the stronger upper bound $\bar{l}\bar{\mathsf{C}}_C < 2(1-1/n)^2 + 2/n^2$ for any connected graph $G$ with $n$ vertices;
see Theorem~\ref{thm:lowerbound}.
We also construct a family of graphs whose values $\bar{l}\bar{\mathsf{C}}_C$ lie arbitrarily close to~$2$;
see Theorem~\ref{thm:Hr21r3}.
We conclude the article by conjecturing that the set of values $\bar{l}\bar{\mathsf{C}}_C$ for all finite connected graphs 
is dense in the interval between~$1$ and~$2$.

\section{Two fundamental bounds and ten lower bounds}
\label{sec:duality}

\noindent
Consider any connected graph $G = (V,E)$ on $n$ vertices and $m$ edges.
For brevity, we will write 
$\bar{\mathsf{C}}_C := \bar{\mathsf{C}}_C(G)$ and 
$\bar{l} := \bar{l}(G)$.

Our first main result are elegantly interlinked lower and upper bounds for $\bar{l}$ and $\bar{\mathsf{C}}_C$.
\begin{thm}\label{thm:CC-l-inequalities}
\[
  1\leq \bar{l}\,\bar{\mathsf{C}}_C < 2\,.
\]
The lower bound is tight and the upper bound is asymptotically tight.
\end{thm}
This theorem will be proved in Section~\ref{sec:proofofmaintheorem}.
The proof shows that $G$ achieves equality $\bar{l}\,\bar{\mathsf{C}}_C = 1$
exactly when the terms $\bar{\mathsf{C}}_C(v)$ are identical for all vertices~$v$;
equivalently, this is when the so-called {\em transmission} $\sum_w d(v,w)$ is independent of the vertex~$v$. 
Graphs with this property were introduced by Handa~\cite{MR1675124} as {\em transmission regular graphs} 
and include cycles, complete graphs, and all other distance-regular graphs 
but which, as shown in~\cite{MR1887380,MR2673007}, are not necessarily regular.
For more information on distance-regular graphs, see, for instance~\cite{MR3684935,MR3640081}.

Theorem~\ref{thm:CC-l-inequalities} provides 
the lower bound $\bar{l}\geq \bar{\mathsf{C}}_C^{-1}$ on the mean distance~$\bar{l}$. 
This complements the many upper bounds on $\bar{l}$
found in the literature; see, for instance~\cite{MR1802606,MR940832,MR485476,MR2886647,MR1105467,MR1664491,MR1931503}.
Conversely, these upper bounds on $\bar{l}$ 
may be combined with the general lower bound $\bar{\mathsf{C}}_C\geq \bar{l}^{-1}$
from Theorem~\ref{thm:CC-l-inequalities} to provide new lower bounds on 
the closeness centrality $\bar{\mathsf{C}}_C$, as described in the following corollary.

\begin{cor}
\label{cor:lower bounds on CC}
Let $G = (V,E)$ be a connected graph on $n$ vertices and $m$ edges with 
minimum vertex degree~$\delta$, maximum vertex degree $\Delta$ and diameter~$\displaystyle D=\max_{u,v\in V}d(u,v)$.
Then
\begin{enumerate}
\setlength{\itemsep}{0mm}
  \item $\displaystyle\bar{\mathsf{C}}_C \geq \Big((n-1)\frac{n-\Delta}{n}\frac{\Delta-1}{\Delta}-\frac{2(m-n+1)}{n(n-1)} + 1\Big)^{-1}$;
  \item $\displaystyle\bar{\mathsf{C}}_C \geq 6n(n-1)\big(3n(n-4)D + 6n(n+2) - 12(m+1) - D^2(D-6) + h(D)\big)^{-1}$;
  \item $\displaystyle\bar{\mathsf{C}}_C \geq (n-1)\bigg(\sum_{k\,:\,Q_k(0)\leq n-1}\Big\lfloor\frac{n(n-1)}{Q_k^2(0)+n-1}\Big\rfloor\bigg)^{-1}$;
  \item $\displaystyle\bar{\mathsf{C}}_C \geq (n-1)n\big(bn(n-1)-2(b-1)m\big)^{-1}$;
  \item $\displaystyle\bar{\mathsf{C}}_C \geq \frac{n-1}{n}\bigg(\Big\lceil\frac{\Delta+\theta_2}{4\theta_2}\ln(n-1)\Big\rceil+\frac{1}{2}\bigg)^{-1}$;
  \item $\displaystyle\bar{\mathsf{C}}_C \geq (n-1)\frac{2(m-n-2)}{\theta_2^{-1}+\cdots+\theta_n^{-1}}$;
  \item $\displaystyle\bar{\mathsf{C}}_C \geq \frac{3}{n+1}$;
  \item $\displaystyle\bar{\mathsf{C}}_C \geq n(n-1)\bigg\lfloor\frac{(n+1)n(n-1)-2m}{\delta+1}\bigg\rfloor^{-1}$;
  \item $\displaystyle\bar{\mathsf{C}}_C \geq \alpha^{-1}$;
  \item $\displaystyle\bar{\mathsf{C}}_C \geq D^{-1}$
\end{enumerate}
where $h(D) = 6-D$ if $n-D$ is odd and $h(D) = 2D$ when $n-D$ is even; 
$\theta_i$ is the $i$th largest eigenvalue of the Laplacian matrix for $G$; 
$b$ is the number of these Laplacian eigenvalues; 
$Q_k$~is the $k$-alternating polynomial defined over them; 
and $\alpha$ is the independence number for~$G$. 
\end{cor}
The first two bounds in Corollary~\ref{cor:lower bounds on CC} follow from upper bounds on $\bar{l}$ proved 
by Gago, Hurajov\'a and Madaras~\cite{MR2886647,MR3125646}; 
the third and fourth bounds follow from upper bounds on $\bar{l}$ by Rodriguez and Yebra~\cite{MR1664491};
and the last six bounds follow from upper bounds on $\bar{l}$ by
Mohar~\cite{MR1105467}, 
Teranishi~\cite{MR1931503}, 
Entringer, Jackson and Snyder~\cite{MR543771},
Beezer, Riegsecker, and Smith~\cite{MR1802606}, 
Chung~\cite{MR940832},
and Doyle and Graver~\cite{MR485476},
respectively. 
The last three of these bounds are tight and are attained exactly when $G$ is the complete graph $K_n$.

\medskip

The lower bound in Theorem~\ref{thm:CC-l-inequalities} was first discovered by the second author in his thesis~\cite{huthesis19}, 
expressed in terms of the {\em betweenness centrality} of $G$. 
This centrality was defined in 1977 by Freeman~\cite{freeman77} as
\[
     \bar{\mathsf{C}}_B
  := \bar{\mathsf{C}}_B(G)
   = \dfrac{1}{n}\sum_{v\in V}\mathsf{C}_B(v)
  \qquad\text{where}\qquad
     \mathsf{C}_B(v) 
   = \sum_{s\neq t}\dfrac{\sigma_{st}(v)}{\sigma_{st}}\,,
\]
where $\sigma_{st}$    is the number of shortest paths between vertices~$s$ and~$t$ and
where $\sigma_{st}(v)$ is the number of such paths that pass through vertex~$v\neq s,t$.
The betweenness centrality $\bar{\mathsf{C}}_B$ and the mean distance $\bar{l}$ 
are related bijectively via the following identity first observed by Gago~\cite{gago06}.
\begin{thm}{\rm\cite{gago06}}\label{thm:Gago1}
$\bar{\mathsf{C}}_B = \frac{1}{2}(n-1)(\bar{l}-1)$.
\end{thm} 
The bound in Theorem~\ref{thm:CC-l-inequalities} thus reflects that 
the betweenness centrality $\bar{\mathsf{C}}_B$ and the closeness centrality $\bar{\mathsf{C}}_C$ 
can be interpreted as dual concepts, 
as noted in~\cite{BrBoFr16}.
The upper bound in Theorem~\ref{thm:CC-l-inequalities} can be used together with known lower bounds on $\bar{l}$
to give upper bounds for $\bar{\mathsf{C}}_C$.
However, we will instead present upper bounds in the following section
that do not rely on Theorem~\ref{thm:CC-l-inequalities}.

\section{Upper bounds on closeness centrality}
\label{sec:upper-bounds}

\noindent 
In Section~\ref{sec:duality}, 
we proved lower bounds on the closeness centrality $\bar{\mathsf{C}}_C(G)$ of a graph~$G=(V,E)$.
In~this section, we prove upper bounds on $\bar{\mathsf{C}}_C(G)$.
To derive the first bound, we make a trivial but useful observation.

\begin{rem}
\label{rem:edge-deletion}
If $e$ is an edge of a connected graph $G$ for which $G-e$ is also connected,
then $\bar{\mathsf{C}}_C(G-e) < \bar{\mathsf{C}}_C(G)$.
\end{rem}

A complete graph is obtained by adding edges to $G$, 
and it is easy to see that complete graphs have closeness centrality $\bar{\mathsf{C}}_C(K_n) = 1$.
Remark~\ref{rem:edge-deletion} therefore implies the following simple upper bound on $\bar{\mathsf{C}}_C(G)$.

\begin{lem}
\label{lem:trivial upper bound}
Let $G$ be a connected graph. 
Then $\bar{\mathsf{C}}_C(G)\leq 1$, and equality is attained exactly when $G$ is a complete graph. 
\end{lem}

We can improve Remark~\ref{rem:edge-deletion} by estimating the reduction of the closeness centrality.
\begin{thm}\label{thm:edge-deletion}
Let $G$ be a connected graph on $n$ vertices and $k$ edges. 
Suppose that $e$ is an edge of $G$ such that $G-e$ is also connected.
Then 
\[
       \bar{\mathsf{C}}_C(G-e) 
  \leq \bar{\mathsf{C}}_C(G)-\frac{2n-2}{n}\Big(\frac{2}{(n-1)(n+2)-2k}-\frac{2}{(n-1)(n+2)+1-2k}\Big)\,.
\]
\end{thm}
\begin{proof}
Let $u$ and $v$ be the end-vertices of $e$.
Then $d_{G-e}(u,v)\geq d_G(u,v)$ for all $u,v\in V$ 
and  $d_{G-e}(x,y) = d_{G-e}(y,x)\geq 2=d_G(x,y)+1 = d_G(y,x)+1$. 
Thus,\vspace*{1mm}
\begin{align*}
  &\quad\;\, \bar{\mathsf{C}}_C(G-e)\,n/(n-1)\\[-.5mm]
  &\leq \frac{1}{1+\displaystyle\sum_{v\in V}d_G(x,v)}
       -\frac{1}{  \displaystyle\sum_{v\in V}d_G(x,v)}
       +\frac{1}{1+\displaystyle\!\sum_{v\in V}d_G(y,v)}
       -\frac{1}{  \displaystyle\sum_{v\in V}d_G(y,v)}
       +\sum_{u\in V}\frac{1}{\displaystyle\sum_{v\in V}d_G(u,v)}\\[-.5mm]
  &=    \bar{\mathsf{C}}_C(G)\frac{n}{n-1}-\Bigg(
          \frac{1}{\displaystyle\sum_{v\in V}d_G(x,v)\Big(1+\!\displaystyle\sum_{v\in V}d_G(x,v)\Big)}
         +\frac{1}{\displaystyle\sum_{v\in V}d_G(y,v)\Big(1+\!\displaystyle\sum_{v\in V}d_G(y,v)\Big)}\Bigg)\,.
\end{align*}
The proof follows from $\displaystyle\sum_{v\in V}d_G(w,v)\leq \frac{1}{2}(n-1)(n+2)-k$ for $w\in\{x,y\}$ \cite[Theorem 2.4]{MR543771}.
\end{proof}
\begin{rem}\label{rem:edge-deletion 2}{\rm 
We can replace $\frac{1}{2}(n-1)(n+2)-k$ with $\displaystyle\max_{w\in V} \sum_{v\in V}d_G(w,v)$ \vspace*{-2mm}
for particular classes of graphs to obtain tighter bounds for those classes.}
\end{rem}

\begin{cor}\label{cor:edge-deletion}
If $G$ is a connected graph on $n$ vertices and $e$ edges, then
\[
  \bar{\mathsf{C}}_C(G)\leq 1-\frac{2}{n}+\frac{4(n-1)}{n\big((n-1)(n+2)-2e\big)}\,.
\]
\end{cor}
\begin{proof}
Obtain $G$ by removing exactly $\binom{n}{2}-e$ edges from $K_n$
and apply Theorem~\ref{thm:edge-deletion}.
\end{proof}

The {\em eccentricity} $\epsilon(v)$ of a vertex $v$ of a graph $G$ is 
the greatest distance between $v$ and any other vertex in~$G$.
The {\em radius} $r(G)$ is the smallest eccentricity $\epsilon(v)$ of a vertex $v$ of~$G$.
For each vertex $v\in V$ and integer $i = 1,\ldots,\epsilon(v)$, 
define $d_i(v) = \big|\{u\in V\::\: d(u,v ) = i\}\big|$.

\begin{thm}
\label{thm:non-trivial upper bound 1}
Let $G=(V,E)$ be a connected graph on $n$ vertices with radius $r \geq 1$. 
Then 
\[  
  \bar{\mathsf{C}}_C(G) \leq \frac{n-1}{n-1 + \binom{r}{2}} \,. 
\]
\end{thm} 

\begin{proof}
If $r=1$, then $\binom{r}{2} = 0$, so we may apply Lemma~\ref{lem:trivial upper bound}.
Suppose then that $r\geq 2$.
Since the identity $d_1(v) + \cdots + d_{\epsilon(v)}(v) = n-1$ is true for each vertex $v\in V$,
the sum
\begin{equation}
\label{equ:sum}
    \sum_{u\in V}  d(u,v)
  = \sum_{i=1}^{\epsilon(v)} id_i(v)
\end{equation} 
is minimal exactly when $d_1(v) = n - \epsilon(v)$ and $d_i(v) = 1$ for each $i = 2, \ldots,\epsilon(v)$.
Therefore, 
\[
       \bar{\mathsf{C}}_C(G) 
  =    \dfrac{1}{n}\sum_{v\in V} \bar{\mathsf{C}}_C(v) 
  \leq \dfrac{1}{n}\sum_{v\in V} \frac{n-1}{n - \epsilon(v) + 2 + \cdots + \epsilon(v)} 
  =    \dfrac{1}{n}\sum_{v\in V} \frac{n-1}{n-1 + \binom{\epsilon(v)}{2}}\,.
\]
Since $\epsilon(v) \geq r$ for each vertex $v\in V$, the theorem follows.
\end{proof}

The next theorem provides a similar bound to that in Theorem~\ref{thm:non-trivial upper bound 1} 
but also takes into consideration the maximal vertex degree~$\Delta$ of~$G$.

\begin{thm}
\label{thm:non-trivial upper bound 2}
Let $G=(V,E)$ be a connected graph on $n$ vertices with 
radius $r \geq 2$ and maximum vertex degree $\Delta$. 
Then 
\[  
  \bar{\mathsf{C}}_C(G) \leq \frac{n-1}{2n-1-\Delta + r(r-3)/2} \,. 
\]
\end{thm} 

\begin{proof}  
Since $d_1(v)\leq \Delta$ for each $v\in V$, 
the sum~(\ref{equ:sum}) is minimal when
$d_1(v) = \Delta$, $d_2(v) = n-1-\Delta-(\epsilon(v)-2)$ and 
$d_i(v) = 1$ for each $3\leq i\leq \epsilon(v)$.
Therefore,
\[
       \bar{\mathsf{C}}_C(G) 
  =    \dfrac{1}{n}\sum_{v\in V} \bar{\mathsf{C}}_C(v) 
  \leq \dfrac{1}{n}\sum_{v\in V} \frac{n-1}{2n-1 - \Delta + \epsilon(v)(\epsilon(v)-3)/2}\,.
\]
Since $\epsilon(v) \geq r$ for each vertex $v\in V$, the theorem follows.
\end{proof}

\begin{rem}{\rm
Note that the lower bounds on $\bar{\mathsf{C}}_C(G)$ in Corollary~\ref{cor:lower bounds on CC}
can be combined with the upper bounds on $\bar{\mathsf{C}}_C(G)$ in
Corollary~\ref{cor:edge-deletion} and 
Theorems~\ref{thm:non-trivial upper bound 1} and~\ref{thm:non-trivial upper bound 2}
to form potentially interesting bounds of the form $A\leq B$
where $A$ is one of the right-hand side expressions in the ten bounds of Corollary~\ref{cor:lower bounds on CC}
and where $B$ is one of the three expressions 
\[
  1-\frac{2}{n}+\frac{4n-4}{n\big((n-1)(n+2)-2e\big)}\,,\qquad
  \dfrac{n-1}{n-1+\binom{r}{2}}\qquad\text{and}\qquad
  \dfrac{n-1}{2n-1-\Delta + r(r-3)/2}\,.    
\]}
\end{rem}

\section{Closeness centralities for particular families of graphs}
\label{sec:families}

\noindent
The previous two sections presented lower and upper bounds on the closeness centrality $\bar{\mathsf{C}}_C(G)$ of a graph~$G$.
This section gives lower and upper bounds on $\bar{\mathsf{C}}_C(G)$ for particular classes of graphs.

Lower bounds on $\bar{\mathsf{C}}_C$ are given for three particular classes of graphs in the following proposition.
The bounds are attained by applying Theorem~\ref{thm:CC-l-inequalities} to upper bounds on $\bar{l}$ proved by Hendry~\cite{MR999903}, Plesnik~\cite{MR732013}, 
and Gago, Hujarov\'a and Madaras~\cite{MR2886647}, respectively.

\begin{thm}
\label{thm:lower bounds on CC for specific graphs}
Let $G = (V,E)$ be a connected graph on $n$ vertices.
Then
\begin{enumerate}
  \item if $G$ is self-complementary, then {\rm(}$n\equiv 0 \text{ or }1\pmod{4}$ and{\rm)}
    \[
      \bar{\mathsf{C}}_C(G) \geq 
        \begin{cases} \displaystyle
          \frac{8n-8}{13n-12} & ,\; n\equiv 0\pmod{4}\\[4mm]\displaystyle
          \frac{8n  }{13n-1 } & ,\; n\equiv 1\pmod{4}\,;
        \end{cases} 
    \]
  \item if $G$ is 2-vertex-connected or 2-edge-connected, 
        then $\displaystyle\bar{\mathsf{C}}_C(G) \geq (n-1)\Big\lfloor\frac{1}{4}n^2\Big\rfloor^{-1}$;
  \item if $T$ is a tree with maximal degree $\Delta$, 
        then $\displaystyle\bar{\mathsf{C}}_C(T) \geq \frac{n\Delta}{2(n-\Delta)(\Delta-1)(n-1)+2}$ .  
\end{enumerate}
The second bound is tight and is achieved when $G$ is the cycle~$C_n$.
\end{thm}

Bounds such as those in Theorem~\ref{thm:lower bounds on CC for specific graphs} and in the preceding sections
are particularly useful if it is not possible to calculate a formula for the closeness centrality.
For some classes of graphs, however, such calculations are possible.
Proposition~\ref{prop:closeness-centrality-for-various-graph-families} below
presents explicit expressions for the closeness centrality $\bar{\mathsf{C}}_C(v)$
for the vertices in twelve particular infinite families of graphs.
The proofs of these expressions are easy and straightforward,
and are therefore omitted.
The expressions and their derivations first appeared in the second and third authors' works~\cite{huthesis19,islamreport20}, 
along with expressions for the closeness and betweenness centralities of similar graphs.
The twelve families of graphs also featured in~\cite{MR3294894} 
where their betweenness centralities were calculated.
Each of the families is well-known, 
except possibly the {\em cocktail party graph} $CP(n)$ and the {\em crown graph} $S^0_n$:
these are the graphs obtained by deleting $n$ vertex-disjoint edges from $K_{2n}$ and $K_{n,n}$, 
respectively.

\def\myscale{.8}

\newcommand\Kn{\begin{tikzpicture}[scale=\myscale]
  \pgfmathsetmacro\nn{5}
  \foreach \na in {1,...,\nn}{\pgfmathsetmacro\nangle{360*(\na-1)/\nn+90}
    \coordinate (\na) at (\nangle:1) {};
    \foreach \nb in {1,...,\na}{\draw (\na) -- (\nb);}}
  \foreach \na in {1,...,\nn}{\draw[fill=lightgray,draw=black] (\na) circle(2pt);}
\end{tikzpicture}}

\newcommand\nKn{\begin{tikzpicture}[scale=\myscale]
  \pgfmathsetmacro\nn{5}
  \foreach \na in {1,...,\nn}{\coordinate (\na) at (360*\na/\nn+306:1);}
  \pgfmathsetmacro\nc{int(\nn-2)}
  \foreach \na in {1,...,\nn}{\foreach \nb in {1,...,\nc}{\draw (\na) -- (\nb);}}
  \foreach \na in {1,...,\nn}{\draw[fill=lightgray,draw=black] (\na) circle(2pt);}
\end{tikzpicture}}

\newcommand\Cn{\begin{tikzpicture}[scale=\myscale]
  \pgfmathsetmacro\nn{5}
  \coordinate (0) at (90:1) {};
  \foreach \na in {1,...,\nn}{\coordinate (\na) at (360*\na/\nn+90:1);
    \pgfmathsetmacro\nb{\na-1}{\draw (\na) -- (\nb);}
  \foreach \na in {1,...,\nn}{\draw[fill=lightgray,draw=black] (\na) circle(2pt);}}
\end{tikzpicture}}

\newcommand\Wn{\begin{tikzpicture}[scale=\myscale]
  \pgfmathsetmacro\nn{5}
  \coordinate (0) at (0,0) {};
  \foreach \na in {1,...,\nn}{\coordinate (\na) at (360*\na/\nn+90:1);
    \pgfmathsetmacro\nb{\na-1}{\draw (0) -- (\na) -- (\nb);}
  \draw (1) -- (\nn);
  \foreach \na in {0,...,\nn}{\draw[fill=lightgray,draw=black] (\na) circle(2pt);}}  
\end{tikzpicture}}

\newcommand\Sn{\begin{tikzpicture}[scale=\myscale]
  \pgfmathsetmacro\nn{5}
  \coordinate (0) at (0,0) {};
  \foreach \na in {1,...,\nn}{\coordinate (\na) at (360*\na/\nn+90:1); \draw (0) -- (\na);}
  \foreach \na in {0,...,\nn}{\draw[fill=lightgray,draw=black] (\na) circle(2pt);}
\end{tikzpicture}}

\newcommand\CPn{\begin{tikzpicture}[scale=\myscale]
  \pgfmathsetmacro\nn{int(6)}
  \foreach \na in {1,...,\nn}{\coordinate (\na) at (360*\na/\nn:1);}
  \foreach \na in {1,...,\nn}{\foreach \nb in {1,...,\na}{\pgfmathsetmacro\csum{int(2*(\na-\nb))}
    \ifnum\csum=\nn{}\else\draw (\na) -- (\nb);\fi;}}
  \foreach \na in {1,...,\nn}{\draw[fill=lightgray,draw=black] (\na) circle(2pt);}
\end{tikzpicture}}

\newcommand\Qthree{\begin{tikzpicture}[scale=\myscale]
  \draw[white] (0,-.195) -- (0,-.195); 
  \foreach \x in {0,1}{
    \foreach \y in {0,1}{
      \foreach \z in {0,1}{
        \coordinate (\x\y\z) at (\x+.71*\y,\z+.71*\y);}}}
  \draw (010) -- (110) -- (111) -- (011) -- cycle;
  \foreach \x in {0,1}{\foreach \z in {0,1}{\draw[fill=lightgray,draw=black] (\x1\z) circle(2pt); \draw (\x0\z) -- (\x1\z);}}
  \draw (000) -- (100) -- (101) -- (001) -- cycle;
  \foreach \x in {0,1}{\foreach \z in {0,1}{\draw[fill=lightgray,draw=black] (\x0\z) circle(2pt);}}
\end{tikzpicture}}

\newcommand\Qthreev{\begin{tikzpicture}[scale=\myscale]
  \draw[white] (0,0) -- (0,0); 
  \foreach \x in {0,1}{
    \foreach \y in {0,1}{
      \foreach \z in {0,1}{
        \pgfmathsetmacro\xx{(\x-.5)*(1+0.71*\y)}
        \pgfmathsetmacro\zz{(\z-.5)*(1+0.71*\y)}
        \coordinate (\x\y\z) at (\xx,\zz)}}}
  \draw (010) -- (110) -- (111) -- (011) -- (010) -- (000) -- (100) -- (101) -- (001) -- (000);
  \draw (100) -- (110) (101) -- (111) (001) -- (011);
  \foreach \x in {0,1}{\foreach \y in {0,1}{\foreach \z in {0,1}{\draw[fill=lightgray,draw=black] (\x\y\z) circle(2pt);}}}
\end{tikzpicture}}

\newcommand\Crn{\begin{tikzpicture}[scale=\myscale]
  \pgfmathsetmacro\nn{4}
  \draw[white] (0,-.5) -- (0,.4); 
  \foreach \na in {1,...,\nn}{\pgfmathsetmacro\x{\na-1-(\nn-1)/2}\foreach \ud in {0,1}{\coordinate (\ud\na) at (\x,\ud);}}
  \foreach \na in {1,...,\nn}{\foreach \nb in {1,...,\nn}{\ifnum\na=\nb{}\else\draw (0\na) -- (1\nb);\fi;}}
  \foreach \na in {1,...,\nn}{\foreach \ud in {0,1}{\draw[fill=lightgray,draw=black] (\ud\na) circle(2pt);}}
\end{tikzpicture}}

\newcommand\bkm{\begin{tikzpicture}[scale=\myscale]
  \pgfmathsetmacro\k{3}
  \pgfmathsetmacro\m{4}
  \draw[white] (0,-.5) -- (0,.4); 
  \foreach \na in {1,...,\k}{\pgfmathsetmacro\x{\na-1-(\k-1)/2}\coordinate (1\na) at (\x,1);}
  \foreach \na in {1,...,\m}{\pgfmathsetmacro\x{\na-1-(\m-1)/2}\coordinate (0\na) at (\x,0);}
  \draw (11)+(0, .05) node()[above] {\footnotesize$u_1$};
  \draw (12)+(0, .05) node()[above] {\footnotesize$\cdots$};
  \draw (13)+(0, .05) node()[above] {\footnotesize$u_k$};
  \draw (01)+(0,-.05) node()[below] {\footnotesize$v_1$};
  \draw (02)+(0,-.05) node()[below] {\footnotesize$v_2$};
  \draw (03)+(0,-.05) node()[below] {\footnotesize$\cdots$};
  \draw (04)+(0,-.05) node()[below] {\footnotesize$v_m$};
  \foreach \na in {1,...,\k}{\foreach \nb in {1,...,\m}{\draw (1\na) -- (0\nb);}}
  \foreach \na in {1,...,\k}{\draw[fill=lightgray,draw=black] (1\na) circle(2pt);}
  \foreach \na in {1,...,\m}{\draw[fill=lightgray,draw=black] (0\na) circle(2pt);}
\end{tikzpicture}}

\newcommand\Pn{\begin{tikzpicture}[scale=\myscale]
  \pgfmathsetmacro\nn{3}
  \pgfmathsetmacro\ni{\nn-2}
  \draw (0,0) -- (\nn-2,0) (\nn-1,0) -- (\nn,0);
  \draw[dashed]  (\nn-2,0) -- (\nn-1,0);
  \foreach \x in {0,...,\nn}{\draw[fill=lightgray,draw=black] (\x,0) circle(2pt);}
  \foreach \x in {0,...,\ni}{\draw (\x,-.33) node(){\footnotesize$v_{\x}$};}
  \draw (\nn-1,-.33) node(){\footnotesize$v_{n-2}$};
  \draw (\nn  ,-.33) node(){\footnotesize$v_{n-1}$};  
\end{tikzpicture}}

\newcommand\Ln{\begin{tikzpicture}[scale=\myscale]
  \pgfmathsetmacro\nn{3}
  \pgfmathsetmacro\ni{\nn-2}
  \foreach \x in {0,...,\nn}{\draw (\x,0) --++ (0,1);}
  \foreach \y in {0,1}{\draw (0,\y) -- (\nn-2,\y) (\nn-1,\y) -- (\nn,\y);
                       \draw[dashed]   (\nn-2,\y) -- (\nn-1,\y);
                       \foreach \x in {0,...,\nn}{\draw[fill=lightgray,draw=black] (\x,\y) circle(2pt);}}
  \foreach \x in {0,...,\ni}{\draw (\x,1.33) node(){\footnotesize$u_\x$};
                             \draw (\x,-.33) node(){\footnotesize$v_\x$};}
   \draw (\nn-1,1.33) node(){\footnotesize$u_{n-2}$};
   \draw (\nn-1,-.33) node(){\footnotesize$v_{n-2}$};
   \draw (\nn  ,1.33) node(){\footnotesize$u_{n-1}$};  
   \draw (\nn  ,-.33) node(){\footnotesize$v_{n-1}$};  
\end{tikzpicture}}

\newcommand\CLn{\begin{tikzpicture}[scale=\myscale]
  \pgfmathsetmacro\nn{5}
  \draw[white] (0,-1) -- (0,-1); 
  \foreach \na in {0,...,\nn}{\coordinate (o\na) at (360*\na/\nn+90:1);
                              \coordinate (i\na) at (360*\na/\nn+90:.61);}
  \foreach \na in {1,...,\nn}{\pgfmathsetmacro\nb{int(\na-1)}\draw (i\na) -- (i\nb) -- (o\nb) -- (o\na);}
  \foreach \na in {1,...,\nn}{\foreach \r in {i,o}{\draw[fill=lightgray,draw=black] (\r\na) circle(2pt);}}
\end{tikzpicture}}

\begin{center}
\begin{figure}[ht]
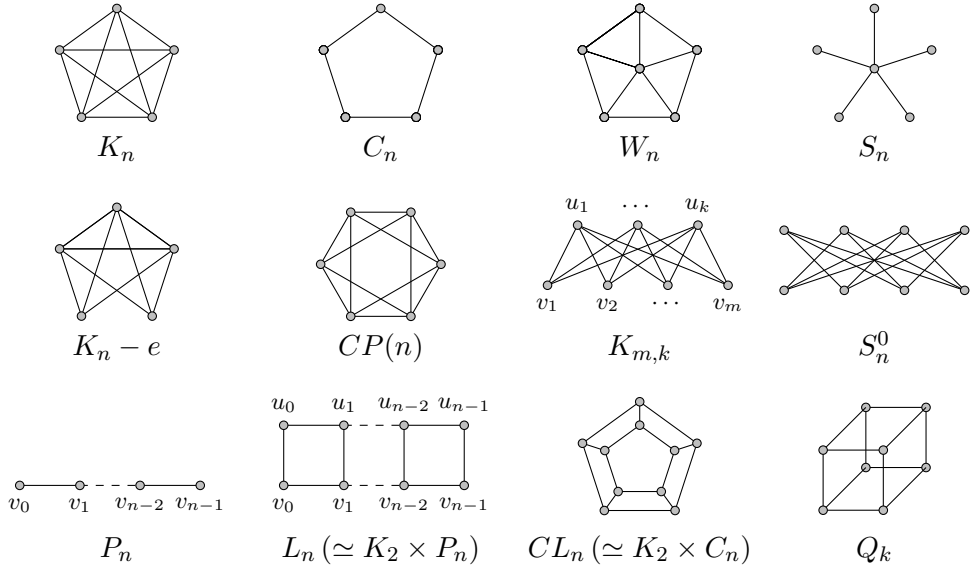

\begin{center}
\begin{tabular}{cccc}
  \Kn     & \Cn     & \Wn       & \Sn      \\
  $K_n$   & $C_n$   & $W_n$     & $S_n$    \\[3mm]
  \nKn    & \CPn    & \bkm      & \Crn     \\
  $K_n-e$ & $CP(n)$ & $K_{m,k}$ & $S^0_n$  \\[3mm]
  \Pn     & \Ln     & \CLn      & \Qthree  \\
  $P_n$   & $L_n\,(\simeq K_2\times P_n)$   & $CL_n\,(\simeq K_2\times C_n)$    & $Q_k$\\
\end{tabular}
\end{center}
\caption{Twelve families of graphs}\label{fig:1}
\end{figure}
\end{center}

\begin{prop}\label{prop:closeness-centrality-for-various-graph-families}
The vertex closeness centralities $\bar{\mathsf{C}}_C(v)$ for all $v\in V$ are given below for the families of graphs $G = (V,E)$ shown 
in Figure~\ref{fig:1}.\\[2mm]
\begin{tabular}{l@{\;}l} 
  Complete graph $K_n${\rm:}          & $\bar{\mathsf{C}}_C(v) = 1$\\
  Cycle    graph $C_n${\rm:}          & $\bar{\mathsf{C}}_C(v) = \frac{n-1}{\lfloor n^2/4\rfloor}$\\[.5mm]
  Wheel    graph $W_n${\rm:}          & $\bar{\mathsf{C}}_C(v) = 1$ if $v$ is the central vertex
                                    and $\frac{n-1}{2n-5}$ otherwise\\[.5mm]
  Star graph $S_n${\rm:}              & $\bar{\mathsf{C}}_C(v) = 1$ if $v$ is the central vertex
                                    and $\frac{n}{2n-1}$ otherwise\\[.5mm]
  Near-complete graph $K_n - e${\rm:} & $\bar{\mathsf{C}}_C(v) = \frac{n-1}{n}$ if $v$ is adjacent to $e$
                                    and $1$ otherwise\\[.5mm]
  Cocktail party graph $CP(n)${\rm:}  & $\bar{\mathsf{C}}_C(v) = \frac{2n-1}{2n}$\\[.5mm]
  Complete bipartite graph $K_{m,k}${\rm:} & $\bar{\mathsf{C}}_C(v) = \begin{cases}
                                       \frac{m+k-1}{m+2k-2} &\text{if } v\in \{u_1,\ldots,u_k\}\\
                                       \frac{m+k-1}{k+2m-2} &\text{if } v\in \{v_1,\ldots,v_m\}
                                     \end{cases}$\\[.5mm]
  Crown graph $S^0_n${\rm:}           & $\bar{\mathsf{C}}_C(v) = \frac{2n-1}{3n}$\\[.5mm]
  Path graph $P_n${\rm:}              & $\bar{\mathsf{C}}_C(v_k) = \frac{4(n-1)}{(2k-n+1)^2 + n^2-1}$ for $k=0,\ldots,n-1$\\[.5mm]
  Ladder graph $L_n${\rm:}            & $\bar{\mathsf{C}}_C(v_k) = \bar{\mathsf{C}}_C(u_k) = \frac{4n-2}{(2k-n+1)^2 + n^2+2n-1}$ for $k=0,\ldots,n-1$\\[.5mm]
  Circular ladder graph $CL_n${\rm:}  & $\bar{\mathsf{C}}_C(v) = \frac{n-1}{2\lfloor n^2/4\rfloor+n}$\\[.5mm]
  Hypercube graph $Q_k${\rm:}         & $\bar{\mathsf{C}}_C(v) = \frac{2^k-1}{k2^{k-1}}$.
\end{tabular}  
\end{prop}

The expressions for $\bar{\mathsf{C}}_C(v)$ in Proposition~\ref{prop:closeness-centrality-for-various-graph-families}
can easily be used to find simple expressions for the closeness centrality $\bar{\mathsf{C}}_C(G)$
of ten of the twelve families of graphs~$G$.

\begin{cor}
\label{cor:closeness-centrality-for-families-of-graphs}
\[
  \begin{array}{l@{\;}l@{\qquad}l@{\;}l}
    \bar{\mathsf{C}}_C(K_n)   &= 1                                & \bar{\mathsf{C}}_C(CP_n)    &= \frac{2n-1}{2n} \\[.5mm]
    \bar{\mathsf{C}}_C(C_n)   &= \frac{n-1}{\lfloor n^2/4\rfloor} & \bar{\mathsf{C}}_C(K_{m,k}) &= \frac{m+k-1}{m+k}\big(\frac{m}{k+2m-2}+\frac{k}{m+2k-2}\big) \\[.5mm]
    \bar{\mathsf{C}}_C(W_n)   &= \frac{n^2-4}{n(2n-5)}            & \bar{\mathsf{C}}_C(S^0_n)   &= \frac{2n-1}{3n} \\[.5mm]
    \bar{\mathsf{C}}_C(S_n)   &= \frac{n^2+2n-1}{2n^2+n-1}        & \bar{\mathsf{C}}_C(CL_n)    &= \frac{n-1}{2\lfloor n^2/4\rfloor+n} \\[.5mm]
    \bar{\mathsf{C}}_C(K_n-e) &= \frac{n^2-2}{n^2}                & \bar{\mathsf{C}}_C(Q_k)     &= \frac{2^k-1}{k2^{k-1}}\,.
  \end{array}
\]
\end{cor}

For the families of paths $P_n$ and ladders $L_n$, 
closed-form expressions for the closeness centrality $\bar{\mathsf{C}}_C(G)$
can be found, for instance by using Maple~\cite{maple2020}.

\begin{cor}
\label{cor:closeness-centrality-for-P_n-and-L_n}
For the graphs $G\in \{P_n,L_n\}$, 
\[
    \bar{\mathsf{C}}_C(G)
  =-\sum_{a,b=\pm1} \frac{abc}{n\hat{n}} \psi\Big(\frac{1 + an + b\hat{n}}{2}\Big)
\]
where $\psi(x)$ is the digamma function and where 
$\hat{n} = \sqrt{1-n^2}$    and $c = \frac{1}{n}$     if $G = P_n$, and 
$\hat{n} = \sqrt{1-2n-n^2}$ and $c = \frac{2n-1}{4n}$ if $G = L_n$.
\end{cor}

Propositions~\ref{prop:cool-lower-and-upper-bounds-on-CC-path} 
and~\ref{prop:cool-lower-and-upper-bounds-on-CC-ladder} 
below provide approximations for $\bar{\mathsf{C}}_C(P_n)$ and $\bar{\mathsf{C}}_C(L_n)$ 
by way of lower and upper bounds which converge asymptotically as~$n$ grows large.

\begin{prop}\label{prop:cool-lower-and-upper-bounds-on-CC-path}
\[
       \frac{4}{n}\sqrt{\frac{n-1}{n+1}}\tan^{-1}\bigg(\sqrt{\dfrac{n-1}{n+1}}\bigg)
  \leq \bar{\mathsf{C}}_C(P_n) 
  \leq \dfrac{\pi}{n}\sqrt{\dfrac{n-1}{n+1}}
     + \frac{n-1}{n\lfloor\frac{n}{2}\rfloor\lceil\frac{n}{2}\rceil}\,.
\]
\end{prop}

\begin{proof}
By Proposition~\ref{prop:closeness-centrality-for-various-graph-families},
\[
    \dfrac{n}{n-1}\,\bar{\mathsf{C}}_C(P_n)
  = \sum_{k=0}^{n-1} f(k) \qquad\text{where}\qquad
    f(x) 
  = \frac{4}{(2x-n+1)^2 + n^2 - 1}\,.
\] 
Since $f$ is non-decreasing for all $x\leq\frac{n-1}{2}$ and is symmetric in $x = \frac{n-1}{2}$\,,
\begin{align*}
       \dfrac{n}{n-1}\bar{\mathsf{C}}_C(P_n) - f\Big(\Big\lfloor{\frac{n}{2}}\Big\rfloor\Big)        
  &=   \sum_{k=0}^{\lfloor{\frac{n}{2}\rfloor-1}}   f(k)
     + \sum_{k=\lfloor{\frac{n}{2}}\rfloor+1}^{n-1} f(k)\\
  &\leq \int_0^{(n-1)/2}f(x)dx+\int_{(n-1)/2}^{n-1}f(x)dx
     =  \int_0^{n-1} f(x)dx\,.
\end{align*}
Therefore by integration, 
\[
    \dfrac{n}{n-1}\bar{\mathsf{C}}_C(P_n)
  - f\Big(\Big\lfloor{\frac{n}{2}}\Big\rfloor\Big) 
  \leq \Biggl[\dfrac{2\tan^{-1}{\biggl(\dfrac{2x-n+1}{\sqrt{n^2-1}}\biggr)}}{\sqrt{n^2-1}}\Biggr]_0^{n-1}
  = \dfrac{4}{\sqrt{n^2-1}}\tan^{-1}\bigg(\dfrac{n-1}{\sqrt{n^2-1}}\bigg)\,.
\]
Hence, 
\[
       \bar{\mathsf{C}}_C(P_n)
    - \frac{n-1}{n}f\Big(\Big\lfloor{\frac{n}{2}}\Big\rfloor\Big) 
  \leq \dfrac{4}{n}\sqrt{\dfrac{n-1}{n+1}}\tan^{-1}\bigg(\sqrt{\dfrac{n-1}{n+1}}\bigg)\,.
\]
Since $\sqrt{\dfrac{n-1}{n+1}}\leq1$, 
\[
       \bar{\mathsf{C}}_C(P_n)
  \leq \dfrac{\pi}{n}\sqrt{\dfrac{n-1}{n+1}}
    + \frac{n-1}{n}f\Big(\Big\lfloor{\frac{n}{2}}\Big\rfloor\Big) 
  =   \dfrac{\pi}{n}\sqrt{\dfrac{n-1}{n+1}}
    + \frac{n-1}{n\lfloor\frac{n}{2}\rfloor\lceil\frac{n}{2}\rceil}\,.
\]
The lower bound is proved similarly.
\end{proof}

The proof of Proposition~\ref{prop:cool-lower-and-upper-bounds-on-CC-ladder}
below is similar 
to that above.

\begin{prop}\label{prop:cool-lower-and-upper-bounds-on-CC-ladder}
\[
  \dfrac{4n-2}{n\sqrt{n^2+2n-1}}\tan^{-1}\biggl(\!\dfrac{n-1}{\sqrt{n^2+2n-1}}\!\biggr)
     + \frac{2}{n^3}\frac{(2n-1)^2}{n^2+2n-1}
  \leq\bar{\mathsf{C}}_C\big(L_n\big)
  \leq\dfrac{\pi}{n}\sqrt{\dfrac{2n-1}{2n+5}}
 +\frac{2}{n(2n+5)}\,.
\]
\end{prop}

\medskip

\noindent
Interestingly, 
Propositions~\ref{prop:cool-lower-and-upper-bounds-on-CC-path} 
         and~\ref{prop:cool-lower-and-upper-bounds-on-CC-ladder}
imply that the closeness centrality measures for paths and ladders
both converge to $\dfrac{\pi}{n}$ as $n$ grows large: 

\medskip

\begin{cor}\label{cor:CCPn-CCLn}
 \[\displaystyle\lim_{n\to\infty} n\bar{\mathsf{C}}_C\big(P_n\big)
  = \lim_{n\to\infty} n\bar{\mathsf{C}}_C\big(L_n\big)
  = \pi\,.
 \]
\end{cor}

\bigskip

\section{Product of closeness centralities and mean distances for graph families}
\label{sec:product-of-closeness-and-mean-distance}

\noindent
By Theorem~\ref{thm:CC-l-inequalities}, 
the product of the closeness centrality and the mean distance of each graph $G$
satisfies the inequality 
$1\leq \bar{l}\,\bar{\mathsf{C}}_C < 2$.
In this section, 
we will calculate the product $\bar{l}\,\bar{\mathsf{C}}_C$ explicitly for specific classes of graphs.
To do so, we first present some easily-calculated mean distances~$\bar{l}(G)$.

\begin{prop}\label{prop:mean-distance-for-various-graph-families}
The mean distances $\bar{l}(G)$ for the families of graphs in Figure~\ref{fig:1} are:
\[
  \begin{array}{l@{\;}l@{\;\;}l@{\;}l@{\;\;}l@{\;}l@{\;\;}l@{\;}l}
    \bar{l}(K_n)   &= 1                                & \bar{l}(S_n)   &= \frac{2n}{n+1}          & \bar{l}(K_{m,k}) &= \frac{m(k+2m-2)+k(m+2k-2)}{(m+k)(m+k-1)} & \bar{l}(L_n)     &= \frac{n+2}{3} \\[.5mm]
    \bar{l}(C_n)   &= \frac{\lfloor n^2/4\rfloor}{n-1} & \bar{l}(K_n-e) &= \frac{n(n-1)+2}{n(n-1)} & \bar{l}(S^0_n)   &= \frac{3n}{2n-1}                          & \bar{l}(CL_n)    &= \frac{2\lfloor n^2/4\rfloor+n}{n-1}\\[.5mm]
    \bar{l}(W_n)   &= \frac{2n-4}{n}                   & \bar{l}(CP_n)  &= \frac{2n}{2n-1}         & \bar{l}(P_n)     &= \frac{n+1}{3}                            & \bar{l}(Q_k)     &= \frac{k2^{k-1}}{2^k-1}\,.
  \end{array}
\]
\end{prop}
For several well-known families of graphs, 
the values of $\bar{l}\,\bar{\mathsf{C}}_C$ attain the lower bound~$1$ in Theorem~\ref{thm:CC-l-inequalities}, 
as the following corollary to Proposition~\ref{prop:mean-distance-for-various-graph-families} shows.
\begin{cor}\label{cor:CC-l-equal to 1}
$\bar{l}\,\bar{\mathsf{C}}_C = 1$ for $G\in\{K_n,C_n,CP_n,S^0_n,CL_n,Q_k\}$.
\end{cor}
For other families, 
$\bar{l}\,\bar{\mathsf{C}}_C$ is always greater than~1 
but converges towards~$1$ as $n$ grows large.
\begin{cor}\label{cor:CC-l-converge to 1}
$\displaystyle\lim_{n\to\infty}\bar{l}\,\bar{\mathsf{C}}_C = 1$ for $G\in\{S_n,K_n-e,W_n\}$.
\end{cor}
Meanwhile, the values of $\bar{l}\,\bar{\mathsf{C}}_C$ converge to unexpectedly interesting constants for 
path and ladder graphs.
\begin{cor}\label{cor:CC-l-converge to pi/3}
$\displaystyle\lim_{n\to\infty}\bar{l}\,\bar{\mathsf{C}}_C = \frac{\pi}{3}$ for $G\in\{P_n,L_n\}$.
\end{cor}
\begin{proof}
Straightforward from Corollary~\ref{cor:CCPn-CCLn} and Proposition~\ref{prop:mean-distance-for-various-graph-families}.
\end{proof}

\begin{cor}\label{cor:CC-l-complete bipartite}
For $G = K_{n-k,k}$, 
\[
  \lim_{n\to\infty}\max_{1\leq k \leq n-1}\bar{l}\,\bar{\mathsf{C}}_C\leq 18-12\sqrt{2}\,.
\]
Indeed,
\[
  \lim_{n\to\infty}\bar{l}\bar{\mathsf{C}}_C = 18-12\sqrt{2}
\]
for 
\[
  k = \Bigg\lfloor\frac{n}{2}-\frac{\sqrt{(3n-4)^2-2(3n-4)\sqrt{2(n-1)(n-2)}}}{2}\Bigg\rfloor\,.
\]
\end{cor}

\begin{proof}
Apply simple calculus and Maple~\cite{maple2020} to Corollary~\ref{cor:closeness-centrality-for-families-of-graphs} 
and Proposition~\ref{prop:mean-distance-for-various-graph-families}. 
\end{proof}

We also consider another family of graphs that provides another unexpectedly interesting asymptotic limit.
For any positive integer $n,k$, 
let $G_{n,k}$ denote the graph 
obtained by attaching a path of $k$ vertices to each vertex of the complete graph $K_n$;
see Figure~\ref{fig:2} for an example. 
\def\myscale{.8}
\newcommand\Gnk{\begin{tikzpicture}[scale=\myscale]
  \pgfmathsetmacro\nn{5}
  \foreach \na in {0,...,\nn}{\coordinate (o\na) at (360*\na/\nn+90:1.4);
                              \coordinate (i\na) at (360*\na/\nn+90:0.7);
                              \coordinate (e\na) at (360*\na/\nn+90:2.1);}                              
  \foreach \na in {1,...,\nn}{\pgfmathsetmacro\nb{int(\na-1)}\draw (i\na) -- (o\na) -- (e\na);}
  \foreach \na in {1,...,\nn}
    {\foreach \nb in {1,...,\na}{\draw (i\na) -- (i\nb);}}
  \foreach \na in {1,...,\nn}{\foreach \r in {i,o,e}{\draw[fill=lightgray,draw=black] (\r\na) circle(2pt);}}
\end{tikzpicture}}
\begin{center}
  \begin{figure}[h!]
    \Gnk
    \caption{$G_{5,2}$}
    \label{fig:2}
  \end{figure}
\end{center}
\begin{prop}\label{prop:lGnk-CCGnk}
\begin{align*}
    \bar{l}\big(G_{n,k}\big)
 &= \frac{(k+1)^2 n-\frac{2}{3}k^2-\frac{4}{3}k-1}{(k+1)n-1}\\\text{and}\quad
    \bar{\mathsf{C}}_C\big(G_{n,k}\big)
 &= \frac{(k+1)n-1}{k+1}\sum_{j=0}^k\frac{2}{j(j+1)+(k-j)(k-j+1)+(n-1)(k+1)(2j+k+2)}\,.
\end{align*}
\end{prop}
\begin{proof}
Let $G_{n,k}$ be written as $\big(\bigcup_{i=1}^n P_{k+1}^i\big)+E$
where $P_{k+1}^i = v_0^i v_1^i\ldots v_k^i$ is a path 
and   $E = \big\{\{v_0^i,v_0^j\}\::\:i,j\in\{1,\ldots,n\},i\neq j\big\}$. 
For any $i\in\{1,\dots,n\}$ and $j\in\{0,\dots,k\}$, 
\begin{align*}
    \sum_{v\in V}d\big(v_j^i,v\big)
  &=\sum_{r=1}^{j-1}d\big(v_j^i,v_r^i\big)+\sum_{r=j+1}^kd\big(v_j^i,v_r^i\big)
   +\sum_{s\neq i}\sum_{r=0}^k d\big(v_j^i,v_r^s\big)\\
  &=\sum_{r=1}^{j-1} r + \sum_{r=j+1}^k(r-j)+(n-1) \sum_{r=0}^k (j+r+1)\\
  &=\frac{j(j+1)}{2}+\frac{(k-j)(k-j+1)}{2}+\frac{(n-1)(k+1)(2j+k+2)}{2}\,.
\end{align*}
The values of $\bar{l}$ and $\bar{\mathsf{C}}_C$ can now be obtained directly from their definition.
\end{proof}
\begin{thm}
Let $G = G_{n,k}$.
For each integer $k\geq 1$, 
$\displaystyle\ln\Big(3-\frac{2}{k+2}\Big)
\leq \lim_{n\to\infty}\bar{l}\,\bar{\mathsf{C}}_C
\leq \ln\Big(3+\frac{2}{k}\Big)$.
Moreover,
\[
  \lim_{k\to\infty}\Big(\lim_{n\to\infty}\bar{l}\,\bar{\mathsf{C}}_C\Big) = \ln 3\,.
\]
\end{thm}
\begin{proof}
It is easy to check that
$\displaystyle\lim_{n\to\infty}\bar{l}\,\bar{\mathsf{C}}_C = \sum_{j=0}^k \frac{1}{2j+k+2}$. 
The bounds follow from the inequalities
$\displaystyle\int_0^{k+1}\hspace*{-4mm}f(x)\,dx \leq \sum_{j=0}^k f(j)\leq \int_0^{k+1}\hspace*{-4mm}f(x-1)\,dx$
where $f(x)=\dfrac{1}{2x+k+2}$.
\end{proof}

\section{Proof of Theorem~\ref{thm:CC-l-inequalities}}
\label{sec:proofofmaintheorem}

\noindent
In order to prove Theorem~\ref{thm:CC-l-inequalities}, 
we first prove the lower bound therein.
\begin{thm}\label{thm:lowerbound}
Let $G$ be a connected graph. 
Then 
\[
  \bar{l}\,\bar{\mathsf{C}}_C\geq 1\,.
\]
\end{thm}
\begin{proof}
Note that
\[
     \sum_{v\in V}\frac{1}{\bar{\mathsf{C}}_C(v)}
   = \frac{1}{n-1}\sum_{v\in V}\sum_{u\in V} d(u,v)
   = n\bar{l}\,,
\]
so, by the Arithmetic-Harmonic Inequality, 
\[
       \bar{l}\,\bar{\mathsf{C}}_C
  =    \bar{l}\,\frac{\displaystyle\sum_{v\in V} \bar{\mathsf{C}}_C(v)}{n} 
  \geq \bar{l}\,\frac{n}{\displaystyle\sum_{v\in V} \frac{1}{\bar{\mathsf{C}}_C(v)}}
  =    1\,.\qedhere
\]
\end{proof}
Next, we prove an upper bound on $\bar{l}\,\bar{\mathsf{C}}_C$ that is stronger than that in Theorem~\ref{thm:CC-l-inequalities}.
\begin{thm}\label{thm:CC-l-inequalities 2}
Let $G$ be a connected graph on $n$ vertices. 
Then 
\[
  \bar{l}\,\bar{\mathsf{C}}_C < 2\Big(1-\frac{1}{n}\Big)^2+\frac{2}{n^2}\,.
\]
\end{thm}
\begin{proof}
Since $d(x,u)+d(u,y)\geq d(x,y)$ for all $x,y,u\in V$,
it holds that, for each $u\in V$, 
\[
  \sum_{x,y\in V-\{u\}\,,\; x\neq y} \hspace*{-5mm} d(x,u)+d(u,y)\;\geq \sum_{x,y\in V-\{u\}}\hspace*{-4mm} d(x,y)\,.
\]
Adding $\sum_{x\in V} d(x,u)+\sum_{y\in V} d(u,y)$ to both sides yields a nice inequality
for each $u\in V$:
\[
  (2n-2)\sum_{v\in V}d(u,v)\geq \sum_{x,y\in V}d(x,y)\,.
\]
Finally, $\bar{l}\,\bar{\mathsf{C}}_C$ is bounded as follows:
\begin{align*}
     \bar{l}\,\bar{\mathsf{C}}_C 
  &= \frac{1}{n^2}\sum_{u\in V}\frac{\displaystyle\sum_{x,y\in V} d(x,y)}{\displaystyle\sum_{v\in V} d(u,v)}
   = \frac{2}{n}
    +\frac{1}{n^2}\sum_{u\in V}\frac{\displaystyle\sum_{x,y\in V-\{u\}}\hspace*{-4mm} d(x,y)}{\displaystyle\sum_{v\in V}d(u,v)}
\leq \frac{2}{n} +\frac{2n-2}{n^2}\frac{\displaystyle\sum_{u\in V}\sum_{x,y\in V-\{u\}}\hspace*{-4mm} d(x,y)}{\displaystyle\sum_{x,y\in V}d(x,y)}\,.
\end{align*}
Since equality is not possible here, 
and since 
\[
  \sum_{u\in V}\sum_{x,y\in V-\{u\}}\hspace*{-4mm} d(x,y) = (n-2)\sum_{x,y\in V} d(x,y)\,,
\]
the proof follows.
\end{proof}

\begin{cor}\label{cor:CC-l-inequality}
Let $G$ be a connected graph. 
Then 
\[
  \bar{l}\,\bar{\mathsf{C}}_C< 2\,.
\]
\end{cor}

Next, we consider one more class of graphs, as follows.
For all positive integers $k,m,n$, 
let $H_{k,m,n}$ be the graph 
obtained by joining edges between all vertices of the complete graph $K_n$ 
to one endpoint of each of $m$ copies of the path $P_{k+1}$; 
see Figure~\ref{fig:3} for an example.
\def\myscale{.8}
\newcommand\Hkmn{\begin{tikzpicture}[scale=\myscale]
  \pgfmathsetmacro\kk{3}
  \pgfmathsetmacro\mm{2}
  \pgfmathsetmacro\nn{5}
  \pgfmathsetmacro\vo{.33}
  \pgfmathsetmacro\vs{2*\vo/(\mm-1)}
  \foreach \na in {0,...,\nn}{\coordinate (k\na) at (360*\na/\nn:1);}
  \foreach \ii in {0,...,\kk}{\foreach \jj in {1,...,\mm}{\coordinate (p\ii\jj) at (2+\ii,\vo+\vs-\jj*\vs);}}
  \foreach \na in {1,...,\nn}{
    \foreach \nb in {1,...,\nn}{\draw (k\na) -- (k\nb);}
    \foreach \jj in {1,...,\mm}{\draw[gray!50!white] (k\na) -- (p0\jj);}}
  \foreach \jj in {1,...,\mm}{\draw (p0\jj) -- (p\kk\jj);}
 \foreach \na in {1,...,\nn}{\draw[fill=lightgray,draw=black] (k\na) circle(3pt);}
  \foreach \ii in {0,...,\kk}{\foreach \jj in {1,...,\mm}{\draw[fill=lightgray,draw=black] (p\ii\jj) circle(3pt);}}
\end{tikzpicture}}
\begin{center}
  \begin{figure}[h!]
    \Hkmn
    \caption{$H_{3,2,5}$}
    \label{fig:3}
 \end{figure}
\end{center}
When $m=1$, this family coincides with the graphs mentioned by Entringer, Jackson and Snyder~\cite[Figure 2.1]{MR543771}. 
Below are semi-explicit expressions for $\bar{l}\big(H_{k,m,n}\big)$ 
and $\bar{\mathsf{C}}_C\big(H_{k,m,n}\big)$; 
we omit the messy but easily-derived explicit forms.
Write $H_{k,m,n}$ as $\big(\bigcup_{i=1}^m P_{k+1}^i\cup K_n\big)+E$, 
where $P_{k+1}^i=v_0^iv_1^i\ldots v_k^i$ is a path 
and   $E = \big\{\{v_0^i,v\}\::\:i\in\{1,\ldots,m\},v\in K_n\big\}$. 
\begin{prop}\label{prop:lHkmn-CCHkmn}
Let $V = V(H_{k,m,n})$.
Then
\begin{align*}
  \bar{l}\big(H_{k,m,n}\big)
  &= \displaystyle\frac{2n(n-1)+nm(k+1)(k+2)+2m\sum_{j=0}^k\sum_{v\in V}d\big(v_j^1,v\big)}
                       {2\big((k+1)m+n\big)\big((k+1)m+n-1\big)}\\
     \bar{\mathsf{C}}_C\big(H_{k,m,n}\big)
  &= \frac{(k+1)m+n-1)}{(k+1)m+n)}\Bigg(\frac{2n}{2n-2+m(k+1)(k+2)}
                                       +\sum_{j=0}^k\frac{m}{\sum_{v\in V}d\big(v_j^1,v\big)}\Bigg)\,.
\end{align*}
\end{prop}

\begin{proof}
For each $u\in V(K_n)$,
\begin{align*}
    \sum_{v\in V}d(u,v)
  = \sum_{v\in V(K_n)}\!\!d(u,v)+\sum_{i=1}^m\sum_{j=0}^k d\big(v,v_j^i\big)
 &= n-1+m\sum_{j=0}^k(j+1)\\
 &= n-1+\frac{m(k+1)(k+2)}{2}\,.
\end{align*} 
For any $i\in\{1,\dots,m\}$ and $j\in\{0,\dots,k\}$, we have
\begin{align*}
     \sum_{v\in V}d\big(v_j^i,v\big)
  &= \sum_{v\in V(K_n)}d\big(v_j^i,v\big)
    +\sum_{r=1}^{j-1}d\big(v_j^i,v_r^i\big)
    +\sum_{r=j+1}^kd\big(v_j^i,v_r^i\big)
    +\sum_{s\neq i}\sum_{r=0}^k d\big(v_j^i,v_r^s\big)\\
  &= \sum_{v\in V(K_n)}(j+1)+\sum_{r=1}^{j-1} r 
    +\sum_{r=j+1}^k(r-j)+(m-1) \sum_{r=0}^k (j+r+2)\\
  &= (n-1)(j+1)+\frac{j(j+1)}{2}+\frac{(k-j)(k-j+1)}{2}+\frac{(m-1)(k+1)(2j+k+4)}{2}\,.
\end{align*}
Since $H_{k,m,n}$ has $(k+1)m+n$ vertices, 
the values of $\bar{l}$ and $\bar{\mathsf{C}}_C$ follow from their definition.
\end{proof}
Consider the examples below, 
calculated by using Maple~\cite{maple2020}:
\[
  H_{10,2,80}           :\; \bar{l}\,\bar{\mathsf{C}}_C      \approx 1.2688\,,\quad
  H_{10^4,2,2\times10^6}:\; \bar{l}\,\bar{\mathsf{C}}_C\approx 1.9547\,,\quad 
  H_{10^6,3,3\times10^9}:\; \bar{l}\,\bar{\mathsf{C}}_C\approx 1.9956\,.
\]
These values suggests that, for some large values of $n$ and $k$, 
the value of $\bar{l}\,\bar{\mathsf{C}}_C$ for some members of this family will get very close to~$2$.
Indeed, this is true, as the following theorem shows.
\begin{thm}\label{thm:Hr21r3}
For $G= H_{r^2,1,r^3}$,
\[
  \displaystyle\lim_{r\to\infty}\bar{l}\,\bar{\mathsf{C}}_C = 2\,.
\]
\end{thm}
\begin{proof}
By Proposition~\ref{prop:lHkmn-CCHkmn}, 
\[
    \bar{l}\big(H_{r^2,1,r^3}\big)
  = \frac{3r^5+4r^4+9r^3+3r^2+3r+2}{3(r^3+r^2+1)(r+1)}
\]
and
\[
    \bar{\mathsf{C}}_C\big(H_{r^2,1,r^3}\big)
  = \frac{r^2(r+1)}{r^3+r^2+1}
    \Bigg(\frac{2r}{r^2+2r+3}+\sum_{j=0}^k\frac{1}{\big(j+\frac{r^3-r^2}{2}\big)^2-\frac{r^2}{4}\big(r^4 - 2r^3 - r^2 - 4r - 2\big)}\Bigg)\,.
\]
Since $\displaystyle\sum_{j=0}^k\frac{1}{\big(j+\frac{r^3-r^2}{2}\big)^2-\frac{r^2}{4}\big(r^4 - 2r^3 - r^2 - 4r - 2\big)}$ is always positive, 
\[
    \bar{l}\,\bar{\mathsf{C}}_C \big(H_{r^2,1,r^3}\big)
  > \bar{l}\big(H_{r^2,1,r^3}\big)\frac{r^2(r+1)}{r^3+r^2+1}\frac{2r}{r^2+2r+3}\,.
\]
The proof now follows since
$\displaystyle\lim_{r\to\infty}\bar{l}\big(H_{r^2,1,r^3}\big)\frac{r^2(r+1)}{r^3+r^2+1}\frac{2r}{r^2+2r+3} = 2$ 
and since  $\bar{l}\,\bar{\mathsf{C}}_C < 2$ by Corollary~\ref{cor:CC-l-inequality}.
\end{proof}

We are now able to prove Theorem~\ref{thm:CC-l-inequalities}.
\begin{proof}[Proof of Theorem~\ref{thm:CC-l-inequalities}]
Apply Theorems~\ref{thm:lowerbound}, \ref{thm:CC-l-inequalities 2} and~\ref{thm:Hr21r3}
and Corollary~\ref{cor:CC-l-equal to 1}.
\end{proof}

Finally, note that Theorem~\ref{thm:Hr21r3} has the following corollary.
\begin{cor}\label{cor:supremum}
For all $\epsilon>0$, 
there exists a connected graph $G$ for which $\bar{l}\,\bar{\mathsf{C}}_C > 2-\epsilon$.
\end{cor}
Partly inspired by this corollary, we suggest the following conjecture.
\begin{conj}
The numbers $\bar{l}\,\bar{\mathsf{C}}_C$ for all connected graphs $G$ 
form a dense subset of the interval~$[1,2)$.
\end{conj}

\section*{ORCID}
\begin{tabular}{ll}
Thomas Britz:           & \url{https://orcid.org/0000-0003-4891-3055}\\
Xin Hu:                 & \url{https://orcid.org/0000-0001-5733-4544}\\
Abdellah Islam:         & \url{https://orcid.org/0000-0002-9659-299X}\\
Hopein Christofen Tang: & \url{https://orcid.org/0000-0003-1707-851X}
\end{tabular}

\section*{Declaration of competing interest}
There is no conflict of interest.

\end{document}